\documentclass[12pt,reqno]{article}

\usepackage[T1]{fontenc}
\usepackage{amsmath}
\usepackage{amsthm}
\usepackage{amscd}
\usepackage{amssymb}
\usepackage{enumerate}
\usepackage{latexsym}
\usepackage{xypic}

\newcommand{\LG}{\mathfrak L}

\newcommand{\dvb}{\|}
\newcommand{\inter}{\overset{\circ}}
\def\epsilon{\varepsilon}
\def\Pr{\mathbb P}
\newcommand{\ssm}{\smallsetminus}

\newcommand{\direction}{\mbox{dir}} 

\newcommand{\im}{\mbox{im}}

\newcommand{\FN}{F_N}   
   
\newcommand{\CVN}{\mbox{CV}_N}   
\newcommand{\barCVN}{\bar{\mbox{CV}}_N}   
\newcommand{\CQ}{{\cal Q}}
\newcommand{\Tobs}{\widehat T^{\mbox{\scriptsize obs}}}
\newcommand{\Tobszero}{\widehat T_{0}^{\mbox{\scriptsize obs}}}
\newcommand{\Tobsone}{\widehat T_{1}^{\mbox{\scriptsize obs}}}

\newcommand{\R}{\mathbb R}

\newcommand{\N}{\mathbb N}
\newcommand{\C}{\mathbb C}
\newcommand{\Hy}{\mathbb H}

\def\strutdepth{\dp\strutbox}
\def \ss{\strut\vadjust{\kern-\strutdepth \sss}}
\def \sss{\vtop to \strutdepth{
\baselineskip\strutdepth\vss\llap{$\diamondsuit\;\;$}\null}}

\def\strutdepth{\dp\strutbox}
\def \sst{\strut\vadjust{\kern-\strutdepth \ssss}}
\def \ssss{\vtop to \strutdepth{
\baselineskip\strutdepth\vss\llap{$\spadesuit\;\;$}\null}}

\def\strutdepth{\dp\strutbox}
\def \ssh{\strut\vadjust{\kern-\strutdepth \sssh}}
\def \sssh{\vtop to \strutdepth{
\baselineskip\strutdepth\vss\llap{$\heartsuit\;\;$}\null}}

\def\qed{\hfill\rlap{$\sqcup$}$\sqcap$\par}
\def\bar{\overline}
\def\tilde{\widetilde}
\def\hat{\widehat}

\newtheorem{thm}{Theorem}[section]
\newtheorem{cor}[thm]{Corollary}
\newtheorem{lem}[thm]{Lemma}
\newtheorem{prop}[thm]{Proposition}

\theoremstyle{definition}
\newtheorem{defn}[thm]{Definition}
\newtheorem*{defn*}{Definition}
\newtheorem{rem}[thm]{Remark}
\newtheorem*{rem*}{Remark}

\theoremstyle{remark}

\numberwithin{equation}{section}

\begin{document}

\title{Non-unique ergodicity, observers' topology and the dual algebraic lamination for $\R$-trees}
%
\author{Thierry Coulbois, Arnaud Hilion, Martin Lustig}

\date{\today }

\maketitle

\begin{abstract}
Let $T$ be an $\R$-tree with a very small action of a free group $F_N$
which has dense orbits.  Such a tree $T$ or its metric completion
$\bar T$ are not locally compact. However, if one adds the Gromov
boundary $\partial T$ to $\bar T$, then there is a coarser {\em
observers' topology} on the union $\bar T \cup \partial T$, and it is
shown here that this union, provided with the observers' topology, is
a compact space $\Tobs$.

To any $\R$-tree $T$ as above a {\em dual lamination} $L^{2}(T)$ has
been associated in \cite{chl1-II}.  Here we prove that, if two such
trees $T_{0}$ and $T_{1}$ have the same dual lamination $L^{2}(T_{0})
= L^{2}(T_{1})$, then with respect to the observers' topology the two
trees have homeomorphic compactifications: $\Tobszero = \Tobsone$.
Furthermore, if both $T_{0}$ and $T_{1}$, say with metrics $d_{0}$ and
$d_{1}$ respectively, are minimal, this homeomorphism restricts to an
$\FN$-equivariant bijection ${T_{0}} \to {T_{1}}$, so that on the
identified set $ {T_{0}} = {T_{1}}$ one obtains a well defined family
of metrics $\lambda d_{1}+(1-\lambda)d_{0}$. We show that for all
$\lambda \in[0,1]$ the resulting metric space $T_{\lambda}$ is an
$\R$-tree.
\end{abstract}


\section*{Introduction}

Geodesic laminations $\LG$ on a hyperbolic surface $S$ are a central
and much studied object in Teichm\"uller theory. A particularily
interesting and sometimes disturbing fact is that there exist minimal
(arational) geodesic laminations that can carry two projectively
distinct transverse measures.  Such minimal {\em non-uniquely ergodic}
laminations were first discovered by W.~Veech \cite{veech} and by
H.~Keynes and D.~Newton \cite{KeyNew}.

\smallskip

The lift $\tilde \LG$ of a geodesic lamination $\LG$, provided with a
transverse measure $\mu$, to the universal covering of the surface $S$
gives rise to a canonical {\em dual} $\R$-tree $T_{\mu}$, which has as
points the leaves of $\tilde \LG$ and their complementary
components. The metric on $T_{\mu}$ is given by the lift of the
transverse measure $\mu$ to $\tilde \LG$, and the action of $\pi_{1}
S$ on $\tilde S \subset \Hy$ induces an action on $T_{\mu}$ by
isometries. We assume that $S$ has at least one boundary component, so
that the fundamental group of $S$ is a free group $\pi_{1} S = \FN$ of
finite rank $N \geq 2$.  For more details see \cite{morg,shalen,chl1-II}.

\smallskip

Distinct measures $\mu$ and $\mu'$ on $\LG$ give rise to dual
$\R$-trees which are not $\FN$-equivariantly isometric. However, the
fact that the two measures are carried by the same geodesic lamination
$\LG$ is sometimes paraphrased by asserting that topologically the two
trees $T_{\mu}$ and $T_{\mu'}$ are ``the same''. We will see below to
what extend such a statement is justified.

\smallskip

In the broader context of very small actions of $\FN$ on $\R$-trees
one can ask whether the analogous phenomenon can occur for an
$\R$-tree $T$ which is not a {\em surface tree}, i.e. it does not
arise from the above construction as dual tree to some measured
lamination on a surface $S$ with $\pi_{1} S \cong \FN$.  An
interesting such example can be found in the Ph.D.-thesis of
M.~Bestvina's student R.~Martin \cite {mart}. Again, in his example
there is a kind of underlying ``geometric lamination'' which is
invariant for different metrics on the ``dual tree'': the novelty in
R.~Martin's example is that the lamination (or rather, in his case,
the {\em foliation},) is given on a finite $2$-complex which is not
homeomorphic to a surface. Such actions have played an important role
in E.~Rips' proof of the Shalen conjecture: they have been termed {\em
geometric} by G.~Levitt, and {\em Levitt} or {\em thin} by others (the
latter terminology being more specific in that it excludes for example
surface trees). They represent, however, by no means the general case
of a very small $\FN$-action on an $\R$-tree,  (see \cite{gl,best}).

\bigskip

In a series of preceding papers \cite{chl1-I, chl1-II, chl1-III} the
tools have been developed to generalize the above two special
situations (``surface'' and ``thin'') sketched above.

\smallskip

As usual, $\partial\FN$ denotes the Gromov boundary of $\FN$ and,
$\partial^2\FN= (\partial\FN)^2\ssm\Delta$ is the double boundary,
where $\Delta$ is the diagonal.  An {\em algebraic lamination} is a
non-empty closed, $\FN$-invariant, flip-invariant subset of
$\partial^2\FN$ (see Definition \ref{lamination}).  This definition
mimicks the set of pairs of endpoints associated to any leaf in the
lift $\tilde \LG$ of a lamination ${\LG} \subset S$ as above.

\smallskip

To any $\R$-tree $T$ with an isometric $\FN$-action we associate a
{\em dual algebraic lamination} $L^{2}(T)$ (without measure), which is
defined as limit of conjugacy classes in $\FN$ with translation length
in $T$ tending to 0.  Moreover, if $T$ is minimal and has dense
orbits, two other equivalent definitions of $L^2(T)$ have been given
in \cite{chl1-II}.

\smallskip

In this paper, we use the definition based on the map $\CQ$ introduced
in \cite{ll4}; $\CQ$ is an $\FN$-equivariant map from $\partial\FN$
onto $\hat T=\bar T\cup \partial T$, where $\bar T$ is the metric
completion of $T$, and $\partial T$ its Gromov boundary.  Details are
given in \S2 below.  The dual algebraic lamination of $T$ is the set
of pairs of distinct boundary points which ``converge'' in $\bar T$ to
the same point $\CQ (X)=\CQ (X')$.  In particular, in the above
discussed case of a surface tree $T = T_{\mu}$, the dual algebraic
lamination $L^{2}(T)$ can be derived directly from the given geodesic
lamination $\LG$, and conversely.

\smallskip

The map $\CQ$ is a strong and useful tool in many circumstances, for
example in the proof of our main result:

\bigskip

\noindent
{\bf Theorem I.}
\label{thm:nonergodicI}
{\em 
Let $T_{0}$ and $T_{1}$ be two $\R$-trees with very small minimal
actions of $\FN$, with dense orbits.  Then the following two
statements are equivalent:

\smallskip
\noindent
(1)  $L^{2}(T_{0}) = L^{2}(T_{1})$.

\smallskip
\noindent
(2) The spaces $\hat T_0$ and $\hat T_1$, both equipped with the
observers' topology, are $\FN$-equivariantly homeomorphic.

Moreover, the homeomorphism of (2) restricts to an $\FN$-equivariant
bijection between $T_0$ and $T_1$.  
}

\bigskip

The observers' topology on the union $\hat T = \bar T \cup\partial T$
is introduced and studied in \S 1 below.  It is weaker than the
topology induced on $\hat T$ by the $\R$-tree's metric, but it agrees
with the latter on segments and also on finite subtrees.  The
difference between the two topologies is best illustrated by
considering an infinite ``multi-pod'' $T^{\infty}$, i.e. a tree which
consists of a central point $Q$ and infinitely many intervals $[P_{i},
Q]$ isometric to $[0, 1] \subset \R$ attached to $Q$. Any sequence of
points $Q_{i} \in [P_{i}, Q)$ converges to $Q$ in the observers'
topology, while in the metric topology one needs to require in
addition that the distance $d(P_{i}, Q)$ tends to 0. In comparison,
recall that in the cellular topology (i.e. $T^\infty$ interpreted as
CW-complex) no such sequence converges.

\bigskip

In the case of a surface $S$ with a marking $\pi_1S=\FN$, there are
several models for the Teichm\"uller space $\mathcal{T}(S)$ and its
Thurston boundary $\partial\mathcal{T}(S)$.  Either it can be viewed
as a subspace of $\Pr\R^{\FN}$, through the lengths of closed
geodesics on the surface $S$, equipped with a hyperbolic structure
(that varies when one moves within $\mathcal{T}(S)$).  In a second
model, the boundary $\partial\mathcal{T}(S)$ can be viewed as the
space of projective measured geodesic laminations. Going from the
second model to the first one is achieved through considering
``degenerated hyperbolic length functions'', each given by integrating
the transverse measure $\mu$ on a geodesic lamination $\LG$ once
around any given closed geodesic.  Alternatively, this amounts to
considering the translation length function of the $\R$-tree $T_\mu$
dual to the measure lamination $(\LG, \mu)$, see \cite{morg}.

\smallskip

We point out that $\partial\mathcal{T}(S)$ is not a convex subset of
$\Pr\R^{\FN}$.  However, the set of projective classes of tranverse
measures on a given minimal (arational) geodesic lamination $\LG$ is a
finite dimensional simplex $\Delta(\LG)$. The extremal points of
$\Delta(\LG)$ are precisely the ergodic measures on $\LG$.

\smallskip

In striking analogy to Teichm\"uller space and its Thurston boundary,
for the free group $\FN$ a ``cousin space'' $\CVN$ has been created by
M.~Culler and K.~Vogtmann \cite{CulVog}.  The points of this {\em
Outer space} $\CVN$ or its boundary $\partial\CVN$ are precisely
given by all non-trivial minimal $\R$-trees, provided with a very
small $\FN$-action by isometries, up to $\FN$-equivariant homothety
(see \cite{cl2}).

\smallskip

Just as described above for $\partial {\cal T}(S)$, there is a
canonical embedding of $\barCVN = \CVN\, \cup \, \partial \CVN$ into
$\Pr\R^{\FN}$, which associates to any homothety class $[T]$ of such
an $\R$-tree $T$ the projective ``vector'' of translation lengths
$\dvb w \dvb_{T}$, for all $w \in \FN$.  (For more detail and
background see \cite{vogt}.)  Hence, for any two homothety classes of
trees $[T_{0}], [T_{1}] \in \barCVN $, there is a line segment
$[T_{0}, T_{1}] \subset \R^{\FN}$ which is given by the set of {\em
convex combinations} of the corresponding translation length
functions.  Again $\barCVN$ is not a convex subspace of $\Pr\R^{\FN}$:
In general, these convex combinations are not length functions that
come from $\R$-trees, and hence the projective image of this line
segment does not lie inside $\barCVN$.  However, we prove:

\bigskip
\noindent
{\bf Theorem II.}
\label{thm:nonergodicII}
{\em Let $T_{0}$ and $T_{1}$ be two minimal $\R$-trees with very small
actions of $\FN$, with dense orbits.  Then statement (1) or (2) of
Theorem I implies:

\smallskip

\noindent
(3) The projectivized image of the segment $[T_{0}, T_{1}] \subset
\R^{\FN}$ of convex combinations of $T_{0}$ and $T_{1}$ is contained
in $\barCVN$.  
}

\bigskip

Our results raise the question of what actually a ``non-uniquely
ergodic'' $\R$-tree is.  Indeed, even this very terminology has to be
seriously questioned.

\smallskip

In the case of trees that are dual to a non-uniquely ergodic surface
lamination, distinct measures on the lamination give rise to
metrically distinct trees.  For the general kind of $\R$-trees $T$
that represent a point in $\partial \CVN$, however, an invariant
measure $\mu$ on the dual algebraic lamination $L^{2}(T)$ (called in
this case {\em a current}, see \cite{kap}) is not directly related to
the metric $d$ on $T$, but much rather defines a dual (pseudo-)metric
$d_{\mu}$ on $T$. It has been shown in \cite{chl1-III} that in general
$d_{\mu}$ is projectively quite different from the original metric
$d$.  Hence we insist on the importance of making a clear distinction
between on the one hand trees with dual algebraic lamination that is
non-uniquely ergodic (in the sense that it supports two projectively
distinct currents) and on the other hand the phenomenon considered in
this paper (see Theorem~I). We suggest the following terminology:

\smallskip

Let $T$ be an $\R$-tree with a minimal $\FN$-action with dense orbits.
$T$ (or rather $\Tobs$) is called {\em non-uniquely ergometric} if there exists a
projectively different $\FN$-invariant metric on $T$ such that the two
observers' topologies coincide.

\bigskip

\noindent{\em Prospective:} 
The work presented in this paper is primarily meant as an answer to a
natural question issuing from our previous work
\cite{chl1-I,chl1-II,chl1-III}, namely: ``To what extend does the dual
algebraic lamination $L^{2}(T)$ determine $T$ ?''  We also hope that
this paper is a starting point for a new conceptual study of
non-unique ergodicity (or rather: ``non-unique ergometricity'') for
$\R$-trees with isometric $\FN$-action.  A first treatment of this
subject, purely in the spirit of property (3) of Theorem II above, has
been given in \S 5 of \cite{guir} (compare also \cite{paul}).  We
believe, however, that there are several additional, rather subtle
topics, which also ought to be adressed in such a study, but which do
not really concern the main purpose of this paper.  To put our paper
in the proper mathematical context, the authors would like to note:

\smallskip
\noindent
(1) The natural question, whether the homeomorphism from part (2) of
Theorem I does extend to an $\FN$-equivariant homeomorphism $T_{0} \to
T_{1}$ with respect to the metric topology, has a negative answer.
Even non-uniquely ergodic surface laminations give already rise to
counterexamples (see \cite{chll}).

\smallskip
\noindent
(2) There are interesting recent results of V.~Guirardel and G.~Levitt
(see \cite{GuiLev}) regarding the converse (under adapted hypotheses)
of the implication given in our Theorem II above.

\smallskip
\noindent
(3) There have been several attempts to introduce ``tree-like
structures'' by purely topological or combinatorial means, which
generalize (or are weaker than) $\R$-trees viewed as topological
spaces. In particular we would like to point the reader's attention to
the work of B.~Bowditch \cite{bow} and that of J.~Mayer, J.~Nikiel and
N.~Oversteegen \cite{mno}, who also consider compactified trees. The
observers' topology seems to be a special case of what they call a
``real tree'', and hence our compactification $\Tobs$ is what they
call a ``dendron''.

\smallskip
\noindent
(4) In the recent book \cite{fj} by C.~Favre and M.~Jonsson one finds
again the observers' topology under the name of ``weak topology'',
introduced for a rather different purpose, namely to study the tree of
valuations for the algebra $\C[[x,y]]$.  Some of the material of our
section 1 can be found in \cite{fj} or already in \cite{mno}, but
translating the references into our terms would be more tedious and
less comfortable for the reader than an independent presentation with
a few short proofs as provided here.

\bigskip

{\em Acknowledgements: This paper originates from a workshop organized
at the CIRM in April 05, and it has greatly benefited from the
discussions started there and continued around the weekly Marseille
seminar ``Teichm\"uller'' (both partially supported by the FRUMAM). We
would in particular like to thank Vincent Guirardel for having pointed
out to us the continuity of the map $\CQ$ with respect to the
observers' topology.}

\bigskip


\section{The observers' topology on an 
$\R$-tree}\label{sec:topobs}

Let $(M, d)$ denote a space $M$ provided with a metric $d$.  The space
$M$ is called {\em geodesic} if any two points $x, y \in M$ are joined
by an arc $[x, y] \subset M$, and this arc is {\em geodesic}: it is
isometric to the interval $[0, d(x, y)] \subset \R$. (Recall that an
arc is a topological space homeomorphic to a closed interval in $\R$,
and an arc {\em joins} points $x$ and $y$ if the homeomorphism takes
the boundary points of the interval to $\{x, y\}$.)

\smallskip

The following remarkable class of metric spaces has been introduced by
M.~Gromov (compare \cite{GhyHar}):

\begin{defn}
\label{def:Gromov}
A metric space $(M, d)$ is called {\em $\delta$-hyperbolic}, with
$\delta \geq 0$, if for any 4 points $x, y, z, w \in M$ one has $(x,
z)_{w} \geq \min\{(x, y)_{w}, (y, z)_{w}\} - \delta$, where $(x,
z)_{w} = \frac{1}{2} (d(w, x) + d(w, z) - d(x, z))$.
\end{defn}

Consider three not necessarily distinct points $P_{1}, P_{2}, P_{3}
\in M$.  We say that $Q \in M$ is a {\em center} of these three points
if for any $i\neq j$ one has $d(P_{i}, P_{j}) = d(P_{i}, Q) + d(P_{j},
Q)$.

\smallskip

\begin{defn}\label{def:Rtree}
An {\em $\R$-tree} $T$ is a metric space which is $0$-hyperbolic and
geodesic.
\end{defn}

Alterna\-tively, a metric space $T$ is an $\R$-tree if and only if any
two points $P, Q \in T$ are joined by a unique arc $[P, Q] \subset T$,
and this arc is geodesic.

\smallskip

We derive directly from the definitions:

\begin{lem}
\label{lem:tripods}
In every $\R$-tree $T$ any triple of points $P, Q, R \in T$ possesses
a unique center $Z \in T$. For any further point $W \in T$ the point
$Z$ is also the center of the triple $W, P, Q$ if and only if one has:

\bigskip

\hfill $(P, Q)_{W} \geq \max\{(P, R)_{W}, (Q, R)_{W}\}$.\qed
\end{lem}

\medskip

For any $\R$-tree $T$ we denote by $\bar T$ the metric completion, by
$\partial T$ the (Gromov) boundary, and by $\hat T$ the union $\bar
T\cup\partial T$.  A point of $\partial T$ is given by a ray $\rho$ in
$T$, i.e. an isometric embedding $\rho: \R_{\geq 0} \to T$. Two rays
$\rho$, $\rho'$ determine the same point $[\rho]$ of $\partial T$ if
and only if their images $\im(\rho)$ and $\im(\rho')$ differ only in a
compact subset of $T$.

\smallskip

The metric on $T$ extends canonically to $\bar T$, and it defines
canonically a topology on $\hat T$ (called below the {\em metric
topology}): A neighborhood basis of a point $[\rho]$ is given by the
set of connected components of $T \ssm \{P\}$ that have non-compact
intersection with $\im(\rho)$, for any point $P \in T$.  We note that
in general $\hat T$ is not compact.

\smallskip

The metric completion $\bar T$ is also an $\R$-tree.  For any two
points $P, Q$ in $\bar T$, the unique closed geodesic arc $[P, Q]$ is
called a {\em segment}.  If $P$ or $Q$ or both are in $\partial T$,
then $[P, Q]$ denotes the (bi)infinite geodesic arc in $\hat T$
joining $P$ to $Q$, including the Gromov boundary point $P$ or $Q$.
    
\smallskip

A point $P$ in $\hat T$ is an {\em extremal point} if $T
\smallsetminus \{P\}$ is connected, or equivalently, if $P$ does not
belongs to the interior $[Q, R] \ssm \{Q, R\}$ of any geodesic segment
$[Q, R]$.  Note that every point of $\partial T$ is extremal, and so
is every point of $\bar T \ssm T$.  We denote by $\inter{T}$ the set
$T$ without its extremal points, and call it the {\em interior tree}
associated to $T$.  Clearly $\inter{T}$ is connected and hence an
$\R$-tree.

\smallskip

For two distinct points $P, Q$ of $\hat T$ we define the {\em
direction $\direction_{P}(Q)$ of $Q$ at $P$} as the connected
component of $\hat T \ssm \{P\}$ which contains $Q$.

\begin{defn}
\label{observerstopology}
On the tree $\hat T$ we define the {\em observers' topology} as the
topology generated (in the sense of a subbasis) by the set of
directions in $\hat T$.  We denote the set $\hat T$ provided with the
observers' topology by $\Tobs$.
\end{defn}

As every direction is an open subset of $\hat T$ (i.e. with respect to
the metric topology), the observers' topology is weaker (= coarser)
than the metric topology.  The identity map $\hat T \to \Tobs$ is
continuous, and isometries of $T$ induce homeomorphisms on $\Tobs$.

\smallskip

The observers' topology has some tricky sides to it which contradict
geometric intuition.  For further reference we note the following
facts which follow directly from the definitions:

\smallskip

\noindent
(a)  An open ball in $\bar T$ is in general not open in $\Tobs$.

\smallskip
\noindent
(b) Every closed ball in $\bar T$ is closed in $\Tobs$.  Note that
closed balls are in general not compact in $\hat T$, but, as will be
shown below, they are compact in $\Tobs$.

\smallskip
\noindent
(c) An infinite sequence of points ``turning around'' a branch point
$P$ (i.e. staying in every direction at $P$ only for a finite time)
converges in $\Tobs$ to $P$.

\medskip

This last property justifies the name of this new topology, which was
suggested by V.~Guirardel: The topology measures only what can be seen
by any set of observers that are placed somewhere in the tree.  We
note as direct consequence of the above definitions:

\begin{rem}
\label{rem:finitetrees}
The restriction of the observers' topology and the restriction of the
metric topology agree on $\partial T$.  Moreover, the two topologies
agree on any finite subtree (=~the convex hull of a finite number of
points) of $\bar T$.
\end{rem}

\begin{lem}
$\Tobs$ is connected and locally arcwise connected.
\end{lem}

\begin{proof} As the observer topology is weaker than the metric topology, any
path for the metric topology is a path for the observers' topology.
As $\hat T$ is arcwise (and locally arcwise) connected, it follows
that $\Tobs$ is arcwise connected, and that elementary open sets
(=~finite intersections of directions) are arcwise connected.
\end{proof}

\begin{prop}
\label{prop:connectedsets}
$\hat T$ and $\Tobs$ have exactly the same connected subsets.  All of
them are arcwise connected for both topologies.
\end{prop}

\begin{proof} 
A connected subset of $\hat T$ is arcwise connected, and therefore it
is also arcwise connected in the observers' topology.

Let $\cal C$ be a connected subset of $\Tobs$, and assume that it is
not connected in the metric topology. Then it is not convex, and hence
there exists points $Q$ and $R$ in $\cal C$ as well as a point $P$ in
$[Q, R]$ which is not contained in $\cal C$.  Now
$U=\mbox{dir}_{P}(Q)$ and $V=\Tobs \smallsetminus
(\mbox{dir}_{P}(Q)\cup \{P\})$ are two disjoint open sets that cover
$\cal C$, with $Q\in U\cap \cal C$ and $R\in V\cap \cal C$. This
contradicts the assumption that $\cal C$ is connected.
\end{proof}

It follows directly from Proposition \ref{prop:connectedsets} that an
extremal point of $\hat T$ is also extremal (in the analogous sense)
in $\Tobs$. In particular we can extend the notion of the {\em
interior tree} to $\Tobs$, and obtain:

\begin{rem}
\label{rem:interiortobs}
The interior trees associated to $\hat T$ and to $\Tobs$ are the same
(as subsets).
\end{rem}

We now observe that in $\inter T$ centers as well as segments have a
very straightforward characterization in terms of directions.

\begin{lem}
\label{lem:center}
(a) A point $Z \in \hat T$ is the center of three not necessarily
distinct points $P_{1}, P_{2}, P_{3} \in \inter T$ if and only if for
any $i \neq j$ the points $P_{i}$ and $P_{j}$ are not contained in the
same connected component of $\hat T \smallsetminus \{Z\}$.

\smallskip
\noindent
(b) A point $R \in \hat T$ belongs to a segment $[P, Q] \subset \inter
T$ if and only if $R$ is the center of the triple $P, Q, R$.  \qed
\end{lem}

The lemma, together with Proposition \ref{prop:connectedsets}, gives
directly:

\begin{prop}
\label{lem:centerpreservation}
\label{lem:segmentspreservation}
Let $T_{0}$ and $T_{1}$ be two $\R$-trees, and assume that there is a
homeomorphism $f: \Tobszero \to\Tobsone$ between the two associated
observers' trees.  Then one has:

\smallskip
\noindent
(a) The center of any three points in ${\inter T}_{0}$ is mapped by
$f$ to the center of the image points in ${\inter T}_{1}$.

\smallskip
\noindent
(b) Any segment $[P, Q]$ in ${\inter T}_{0}$ is mapped to the segment
$[f(P), f(Q)]$ in ${\inter T}_{1}$.  \qed
\end{prop}

Let $(P_n)_{n\in\N}$ be a sequence of points in $\hat T$, and for some
{\em base point} $Q \in \hat T$ consider the set $I_{m} =
\underset{n\geq m}{\bigcap} \, [Q, P_n] \subset\hat T$.  We note that
$I_{m}$ is a segment $I_{m} = [Q, R_{m}]$ for some point $R_{m} \in
\hat T$, and that $I_{m} \subset I_{m+1}$ for all $m \in \N$. Hence
there is a well defined limit point $P = \underset{m \to \infty}{\lim}
\, R_{m}$ (with respect to the metric topology, and thus as well with
respect to the weaker observers' topology), called the {\em inferior
limit from $Q$ of the sequence $(P_{n})_{n \in \N}$}, which we denote
by $P = \underset{n \to \infty}{\liminf_{Q}} \, P_{n}$.
Alternatively, $P$ is characterized by:
$$[Q, P]=\overline{\bigcup_{m=0}^\infty\bigcap_{n\geq m}[Q, P_n]}.$$

It is important to notice that, without further restrictions, the
inferior limit $P$ from some point $Q$ of the sequence
$(P_n)_{n\in\N}$ is always contained in the closure of the convex hull
of the $P_{n}$, but its precise location does in fact depend on the
choice of the base point $Q$.  However, one obtains directly from the
definition:

\begin{lem}
\label{lem:weakconvergence}
Let $(P_{k})_{k \in \N}$ be a sequence of points on $\Tobs$, and let
$D$ be any direction of $\hat T$. Then one has:

\smallskip
\noindent
(a) If all $P_{k}$ are contained in $D$, then for any $Q \in \Tobs$
the inferior limit $\underset{k \to \infty}{\liminf_{Q}} \, P_{k}$ is
contained in the closure $\bar D$ of $D$.

\smallskip
\noindent
(b) If for some $Q \in \Tobs$ the limit inferior $\underset{k \to
\infty}{\liminf_{Q}} \, P_{k}$ is contained in $D$, then infinitely
many of the $P_{k}$ are contained in $D$ as well.

\smallskip
\noindent
(c) If $\underset{k \to \infty}{\liminf_{Q}} \, P_{k}$ lies in $D$ and
if the point $Q$ is not contained in $D$, then all of the $P_{k}$ will
eventually be contained in $D$ as well.  \qed
\end{lem}

\begin{lem}
\label{lem:convergence}
If a sequence of points $P_{n}$ converges in $\Tobs$ to some limit
point $P \in \Tobs$, then for any $Q \in \hat T$ one has:
$$P = \underset{n \to \infty}{\liminf}{}_{Q} \,\, P_{n}$$
\end{lem}

\begin{proof}
From the definition of the topology of $\Tobs$ it follows that any
direction $D$ in $\Tobs$ that contains the limit $P$ will contain all
of the $P_{n}$ with $n$ sufficiently large.  From Lemma
\ref{lem:weakconvergence} (a) it follows that for any $Q \in \Tobs$
the point $R = \underset{n \to \infty}{\liminf}{}_{Q} \,\, P_{n}$ is
contained in the closure $\bar D$, which proves the claim.
\end{proof}

\medskip

We conclude this section with the following observation, which will be
used in section 2, but may also be of independent interest.  Note
that, since any metric space which contains a countable dense subset
is separable, any $\R$-tree $T$ with an action of a finitely generated
group by isometries is separable, if $T$ is minimal or has dense
orbits (see Remark \ref{minimalactions}).

\begin{prop}
\label{prop:Hausdorff}
$\Tobs$ is Hausdorff.  Moreover, if $T$ is separable, then $\Tobs$ is
separable and compact.
\end{prop}

\begin{proof} 
It follows directly from the definition that $\Tobs$ is Hausdorff.
Assume now that $T$ is separable.  It thus contains a countable dense
subset $\chi_{0}$, and hence also a countable subset $\chi$ with the
property that $\chi$ intersects all non-trivial geodesics of $\hat T$:
Such a $\chi$ is given for example as the set of midpoints of any pair
from $\chi_{0} \times \chi_{0}$.

We consider the set of all directions of the form $\direction_{P}(Q)$
with $P, Q \in \chi$, and their finite intersections. It is not hard
to see that this is a countable set which is an open neighborhood
basis for the topology of $\Tobs$.

We now prove that in this case $\Tobs$ is compact.  Let
$(P_n)_{n\in\N}$ be a sequence of points in $\Tobs$ and let
$(D_i)_{i\in\N}$ be a countable family of directions that generates
the open sets of $\Tobs$.  By extracting a subsequence of
$(P_n)_{n\in\N}$ we can assume that for each direction $D_i$ the
sequence $(P_n)_{n\in\N}$ is eventually inside or outside of $D_i$.
We now fix some point $Q \in \Tobs$ and consider the limit inferior $P
= \underset{n \to \infty}{\liminf}{}_{Q} \, P_{n}$ from $Q$. It
follows from Lemma \ref{lem:weakconvergence} (b) that every direction
$D_{i}$ that contains $P$ must contain infinitely many of the $P_{n}$,
and hence, by our above extraction, all but finitely many of them.
This means that the sequence $(P_n)_{n\in\N}$ converges in $\Tobs$ to
$P$.
\end{proof}


\section{The map $\CQ$ and the observers' topology}\label{sec:Q}

From now on let $T$ be an $\R$-tree with a very small minimal action
of a free group $F_N$ by isometries, and assume that some (and hence
any) $\FN$-orbit of points is dense in $T$.

\begin{rem}\label{minimalactions}
An (action on an) $\R$-tree is {\em minimal} if there is no proper invariant subtree.
The ``minimal'' hypothesis is very natural as every $\R$-tree in $\bar{\mbox{CV}}_{N}$ is minimal. 
A minimal tree $T$ is equal to its interior $\inter T$.
Note also that the interior of a tree with dense orbits is minimal. 
\end{rem}

For such trees there is a canonical map $\CQ: \partial F_N\to\hat T$
which has been defined in several equivalent ways in \cite{ll4,ll3}.
Here we use the following definition, which emphasizes the link with
the observers' topology.

\begin{rem}[\cite{ll4,ll3}]\label{rem:Qliminf}  
For all $X$ in $\partial F_N$, for any sequence $(w_i)_{i \in \N}$ in
$\FN$ which converges to $X$ and for any point $P$ of $T$, the point
$$\CQ (X) = \underset{i\to \infty} {\liminf}{}_{P} \, \, w_{i} P \in\hat T$$ 
is independent from the choice of the sequence
$(w_i)_{i\in\N}$ and from that of the point $P$.
\end{rem}

To get some intuition and familiarity with the map $\CQ$ the reader is
refered to \cite{chl1-II}.  The following fact was pointed out to us
by V.~Guirardel:

\begin{prop}
The map $\CQ: \partial \FN \to \Tobs$ is continuous.
\end{prop}

\begin{proof}
We consider any family of elements $X_{k} \in \partial \FN$ that
converges to some $X \in \partial \FN$, with the property that the
sequence of images $\CQ(X_{k})$ converges in $\Tobs$ to some point $Q
\in \Tobs$.  Since $\partial \FN$ and $\Tobs$ are compact (see
Proposition \ref{prop:Hausdorff}), it suffices to show that for any
such family one has $Q = \CQ(X)$. We suppose this is false, and
consider a point $S$ in the interior of the segment $[Q, \CQ(X)]$.

We then consider, for each of the $X_{k}$, a sequence of elements
$w_{k,j} \in \FN$ that converges (for $j \to \infty$) to $X_{k}$. It
follows from the definition of $Q$ that for large $k$ the point
$\CQ(X_{k})$ must be contained in $D = \direction_{S} (Q)$. But then,
by Remark \ref{rem:Qliminf} and Lemma \ref{lem:weakconvergence} (c),
for large $j$ and any $P$ outside $D$ the point $w_{k,j} P$ must also
be contained in $D$.  Hence there exists a diagonal sequence $w_{k,
j(k)}$ which converges to $X$ where all $w_{k,j(k)} P$ are contained
in $D$.  But then Remark \ref{rem:Qliminf} and Lemma
\ref{lem:weakconvergence} (a) implies that $\CQ(X)$ is contained in
$\bar D$, a contradiction.
\end{proof}

Clearly, the map $\CQ$ is $\FN$-equivariant.  Moreover, for the
convenience of the reader, we include a (new) proof of the following
result.

\begin{prop}[\cite{ll4}]
The map $\CQ: \partial \FN \to \Tobs$ is surjective.
\end{prop}

\begin{proof}
By the previous proposition, the image of $\CQ$ is a compact
$\FN$-invariant subset of $\Tobs$.  By hypothesis, $\FN$-orbits are
dense in $T$ for the metric topology.  This implies that $\FN$-orbits
are dense in the metric completion $\bar T$.  Therefore, $\FN$-orbits
are dense in $\bar T$ for the weaker observers' topology and the
$\FN$-orbit of any point in $\bar T$ is dense in $\Tobs$.  It only
remains to prove that the image of $\CQ$ contains a point in $\bar T$.
This is an easy consequence of the fact that the action of $\FN$ is
not discrete.
\end{proof}

\smallskip

Consider the {\em double boundary}
$\partial^2\FN=\partial\FN\times\partial\FN \ssm\Delta$, where
$\Delta$ stands for the diagonal. It inherits canonically from
$\partial \FN$ a topology, an action of $\FN$, and an involution
(called {\em flip}) which exchanges the left and the right factor. For
more details see \cite{chl1-I}, where the following objects have been
defined and studied.

\begin{defn}
\label{lamination}
A non-empty subset $L^{2}$ of $\partial^2\FN$ is an {\em algebraic
lamination} if its closed, $\FN$-invariant and invariant under the
flip involution.
\end{defn}

\smallskip

In \cite{chl1-II} for any 
$\R$-tree $T$ with isometric $\FN$-action the {\em dual algebraic 
lamination} $L^{2}(T)$ has been defined as the set of all 
accumulation points of any family of conjugacy 
classes with translation length on $T$ that tends to $0$.
If the $\FN$-orbits are dense in $T$, then it is proven in 
\cite{chl1-II} that $L^{2}(T)$ is given alternatively by:
$$L^2(T)=\{ (X,X')\in\partial^2 F_N\,|\, \CQ(X)=\CQ(X')\}$$ Here we
focus on the equivalence relation on $\partial F_N$ whose classes are
fibers of $\CQ$, and we denote by $\partial F_N/L^2(T)$ the quotient
set.  The quotient topology on $\partial F_N/L^2(T)$ is the finest
topology such that the natural projection $\pi:\partial F_N\to\partial
F_N/L^2(T)$ is continuous.  The map $\CQ$ splits over $\pi$, thus
inducing a map $\varphi: \partial F_N/L^2(T) \to \Tobs$ with $\CQ =
\varphi \circ \pi$, as represented in the following diagram:
$$\xymatrix{
    \partial F_N \ar@{->>}[rd]^\pi \ar@{->>}[dd]^\CQ \\
    & \partial F_N/L^2(T) \ar[ld]^\varphi_\cong\\
    \Tobs 
  }
$$

\smallskip

By definition of $\partial F_N/L^2(T)$ and by the surjectivity of
$\CQ$, the map $\varphi$ is a bijection. The maps $\CQ$ and $\pi$ are
continuous, and, by virtue of the quotient topology, so is $\varphi$.
As $\Tobs$ is Hausdorff and $\varphi$ is continuous, $\partial
F_N/L^2(T)$ must also be Hausdorff.  Since $\pi$ is onto and
continuous, and $\partial F_N$ is compact, it follows that $\partial
F_N/L^2(T)$ is compact.  Now $\varphi$ is a continuous surjective map
whose domain is a compact Hausdorff space, which shows:

\begin{cor}  
\label{cor:homeo}
The map $\varphi: \partial F_N/L^2(T) \to \Tobs$ is a homeomorphism.
\qed
\end{cor}

This shows that $\Tobs$ is completely determined by the dual algebraic
lamination $L^{2}(T)$ of the $\R$-tree $T$.  As the above defined maps
$\pi, \varphi$ and $\CQ$ are all $\FN$-equivariant, we obtain:

\begin{prop}\label{prop:l2tobs}
Let $T_{0}$ and $T_{1}$ be two $\R$-trees with very small actions of
$F_N$, with dense orbits.  If $L^2(T_{0})=L^2(T_{1})$, then
$\Tobszero$ and $\Tobsone$ are $\FN$-equivariantly homeomorphic.  This
homeomorphism commutes with the canonical maps $\CQ_{0}: \partial \FN
\to \Tobszero$ and $\CQ_{1}: \partial \FN \to \Tobsone$:
$$\xymatrix{
&\partial F_N/L^2(T) \ar[ld]_\cong \ar[rd]^\cong\\
\Tobszero&&\Tobsone
}$$
\qed
\end{prop}

\noindent{\em Proof of Theorem~I.}  The statement of
Proposition~\ref{prop:l2tobs} gives directly the implication from (1)
to (2) in Theorem~I of the Introduction, and, in fact, it seems
slightly stronger.  However, it follows directly from the definition
of the map $\CQ$ and Remark~\ref{rem:Qliminf} that any
$\FN$-equivariant homeomorphism as in Proposition~\ref{prop:l2tobs}
also satisfies the corresponding commutative diagram. In particular,
the converse implication from (2) to (1) in Theorem~I is then a direct
consequence of Corollary \ref{cor:homeo}.

The last part of Theorem~I follows from the fact that $T_0$ and $T_1$
are equal to their interior (see Remark \ref{minimalactions}).  \qed


\section{The proof of   Theorem II}\label{sec:proof}

Let $T_{0}$ and $T_{1}$ be two $\R$-trees with very small
$\FN$-actions with dense orbits, and assume that
$L^2(T_{0})=L^2(T_{1})$.  Then by Proposition \ref{prop:l2tobs} the
associated observers' trees $\Tobszero$ and $\Tobsone$ are
$\FN$-equivariantly homeomorphic.

Through the homeomorphism we identify $\Tobszero$ and $\Tobsone$.
This set is equipped with three topologies (the two metric topologies
and the observers' topology).  In section \ref{sec:topobs} we have
proved that they all have the same connected subsets.  In particular,
they have the same interior tree $\inter{T}$ (see Remark
\ref{rem:interiortobs}), which is also the interior tree of $T_{0}$
and of $T_{1}$.  On this interior tree $\inter T$ both, the metric
$d_{0}$ from $T_{0}$ and the metric $d_{1}$ from $T_{1}$, are well
defined (and finite).

Since any non-negative linear combination of two metrics on the same
space defines a new metric on this space, we can define, for any
$\lambda$ in $[0,1]$, the distance
$$d_\lambda=\lambda d_{1} + (1-\lambda)d_{0}$$ 
on $\inter{T}$, to obtain a metric space ${\inter T}_{\lambda}$.  It
is immediate that $F_N$ acts on ${\inter T}_{\lambda}$ by isometries.

\begin{prop}
\label{linearcombination}
For any $\lambda \in [0,1]$ the metric space ${\inter T}_{\lambda}$ is
a an $\R$-tree.
\end{prop}

\begin{proof}
By Proposition \ref{lem:segmentspreservation} (a) the center of any
triple of points with respect to $d_{0}$ is also the center with
respect to $d_{1}$.  By Lemma \ref{lem:tripods} the three Gromov
products of any triple of points $P, Q, R \in \inter T$ with respect
to a fourth point $W \in \inter T$ are either all three equal for
both, $d_{0}$ and $d_{1}$, or else the maximal one comes from the same
pair for both metrics, and hence the other two pairs have identical
Gromov product with respect $d_{0}$ and $d_{1}$, by Definition
\ref{def:Gromov}.  In both cases the inequality from Definition
\ref{def:Gromov} follows directly for $d_{\lambda}$, so that ${\inter
T}_{\lambda}$ is $0$-hyperbolic.  Furthermore, Proposition
\ref{lem:segmentspreservation} (b) assures us that in ${\inter T}$,
and hence in any ${\inter T}_{\lambda}$, for any two points $P, Q \in
{\inter T}$ there is a well defined segment $[P, Q]$ which agrees with
the segment coming from $T_{0}$ as well with that from $T_{1}$.  By
Remark \ref{rem:finitetrees} the topology on such a segment is the
same for the three topologies carried by $\inter T$, and hence it also
agrees with the topology given by any of the $T_{\lambda}$.  Thus the
$d_{\lambda}$-metric gives an isometry of this segment to the interval
$[0, d_{\lambda}(P, Q)] \subset \R$.  This shows that ${\inter
T}_{\lambda}$ is a $0$-hyperbolic geodesic space, i.e. an $\R$-tree.
\end{proof}

\noindent{\em Proof of Theorem~II.}
Notice first that the assumption in Theorem~II, that both trees 
$T_{0}$ and $T_{1}$ are minimal, implies that both agree with their 
interior subtree (compare Remark \ref{minimalactions}), and hence 
both can be identified canonically with $\inter T$ as above. We now 
apply Proposition \ref{linearcombination} and observe that a linear 
combination of the two metrics $d_{0}$ and $d_{1}$ on $\inter T$ 
implies directly that the corresponding translation length functions 
are given by the analogous linear combination.  This establishes 
statement (3) from Theorem~II as a direct consequence.
\qed


\bibliographystyle{alpha}


\noindent\texttt{thierry.coulbois\string @univ-cezanne.fr\\
arnaud.hilion\string @univ-cezanne.fr\\
martin.lustig\string @univ-cezanne.fr\\
}
Math\'ematiques -- LATP\\
Universit\'e Paul C\'ezanne -- Aix-Marseille III\\
av. escadrille Normandie-Ni\'emen\\
13397 Marseille 20\\ 
France

\end{document}